# KLEIN'S DOUBLE DISCONTINUITY REVISITED:
# WHAT USE IS UNIVERSITY MATHEMATICS
# TO HIGH SCHOOL CALCULUS?

**Abstract** – Much effort and research has been invested into understanding and bridging the 'gaps' which many students experience in terms of contents and expectations as they begin university studies with a heavy component of mathematics, typically in the form of calculus courses. We have several studies of bridging measures, success rates and many other aspects of these "entrance transition" problems. In this paper, we consider the inverse transition, experienced by university students as they revisit core parts of high school mathematics (in particular, calculus) after completing the mandatory undergraduate mathematics courses. To what extent does the "advanced" experience enable them to approach the high school calculus in a deeper and more autonomous way? To what extent can capstone courses support such an approach? How could it be hindered by deficiencies in the students' "advanced" experience? We present a theoretical framework for an analysis of these questions, as well as a number of critical observations and reflections on how they appear in our own institutional context.

# 1. A CLASSICAL PROBLEM REVISITED

In an often-quoted preface, Klein (1908/1932, p. 1) observed that students face a "double discontinuity" as they move from high school to university, then back again to a career as schoolteachers:

> The young university student found himself, at the outset, confronted with problems, which did not suggest, in any particular, the things with which he had been concerned at school. Naturally he forgot these things quickly and thoroughly. When, after finishing his course of study, he became a teacher, he suddenly found himself expected to teach the traditional elementary mathematics in the old pedantic way; and, since he was scarcely able, unaided, to discern any connection between this task and his university mathematics, he soon fell in with the time honoured way of teaching, and his university studies remained only a more or less pleasant memory which had no influence upon his teaching.

The first 'discontinuity' concerns the well-known problems of transition which students face as they enter university, a main theme in research on university mathematics education (see e.g. Gueudet, 2008). The second 'discontinuity' concerns those who which return to school as teachers and the (difficult) transfer of *academic knowledge gained at university* to *relevant knowledge for a teacher*.

Since Klein's days, and particularly in the past few decades, much research has been devoted to mathematics teacher knowledge, and not least on the contributions of initial teacher education (see e.g. Ball, Hill & Bass, 2005; Evens & Ball 2009; Buchholtz et al., 2013). The vast majority of these studies focus on the teachers' preparation and knowledge in primary and lower secondary level mathematics, where the distance to university level mathematics is evident. In this paper our focus will be on future teachers of calculus at the high school (upper secondary) level, and in particular on how an academically oriented bachelor program in pure mathematics may (or may not) contribute knowledge which is relevant to the task of teaching high school calculus. At this level, we have a smaller gap between university mathematics and the mathematics to be taught by the teacher. Indeed it is common for the initial education of teachers at this level to include a substantial amount of pure mathematics courses from a university, which is only natural if we agree with Klein that

> …the teacher's knowledge should be far greater than that which he presents to his pupils. He must be familiar with the cliffs and the whirlpools in order to guide his pupils safely past them. *(Klein, 1908/1932, p. 192)*

But even in this case, there is no reason to expect an automatic transfer (or "trickle-down theory", in the terms of Wu, 2011) between the advanced knowledge gained at university and the tasks of teaching calculus in a high school.

Here are some of the reasons why the preparation of teachers is particularly interesting to us in this particular case:

- In high school, calculus is one of the most advanced topics, and it is usually taught in a quite informal style, leaving the teacher with delicate choices and tasks of explanation (for example in relation to notions of limit, as studied by Barbé et al., 2005);
- The "cliffs and whirlpools" of this topic can indeed be usefully approached using elements of typical undergraduate courses in pure mathematics, particularly in real analysis and algebra, but it could still be necessary to learn such an approach explicitly at university;
- Calculus is particularly affected by the increasing use of symbolic calculation devices in high schools in many countries, and this use in itself leads to important challenges for teaching (see e.g. Guin et al., 2005).

As a result, Klein's problem has significant and critical aspects for the specific case of high school teachers and calculus in our time.

Klein's own *answer* to the problem was that indeed, *university instruction* [must take into account] *the needs of the school teacher* (ibid., p. 1) His proposition for so doing was a series of lectures specifically designed to *help them see the mutual connection between problems in the various fields (…) and more especially to emphasize the relation of these problems to those of school mathematics* (ibid., p. 1-2). He emphasizes that to follow these lectures, students should already be *acquainted with the main features of the chief fields of mathematics* (ibid. p. 1), as a result of their previous university studies. In essence, he thus advocates the introduction of what is known today as *capstone courses* for mathematics teachers (we return to this notion in Section 2).

Klein's proposition implies two concrete problems: how could such a course detect and remedy gaps between students' knowledge and relevant knowledge for the teacher, and what is the appropriate 'higher standpoint' required before the capstone course?

To approach and sharpen these questions, we present (in Section 2) a reformulation of Klein's "double discontinuity", based on the anthropological theory from the didactics of mathematics. In Section 3, we explain the context in which we have studied the questions, namely that of a capstone course within a mainstream university programme in pure mathematics; we also present our methodology, based on the theoretical model, for analysing students' performance on selected tasks from the course. This leads to identifying a number of principal challenges which arise for such a course in the setting of a contemporary academic mathematics programme and which are presented in Section 4, with some wider implications and perspectives for research being exposed in Section 5.

Finally, we offer some remarks about how our paper is situated in, and contributes to, the existing literature. The topic is certainly related to both mathematics teacher education and secondary level mathematics education, but our study focuses on a specific problem for university mathematics education, namely that of using students' "advanced" experience to gain a deeper insight into calculus at high school level. Our paper contributes to the existing literature on this problem by providing a modern theoretical model for the transition to be achieved in students' relation to mathematics, and by showing how the model can be used to analyse students' work in a capstone course where the "advanced" experience is an bachelor programme in pure mathematics (as found in many universities). The fact that such an experience does not automatically ensure "deep" knowledge of elementary mathematics is not new. For instance, Buchholtz et al. (2013) presented a systematic study of student knowledge of similar kind in Germany, China, Hong Kong and South Korea, using a (mainly quantitative) test design. That study, based on a diagnostic test, gives a global picture of the shortcomings, which motivate capstone courses. Our study aims to identify specific challenges which students meet in more complex tasks from a course set up specifically to explore high school calculus while drawing on their "advanced" experience.

2. THEORETICAL FRAMEWORK

We find it useful to relate Klein's "double discontinuity" to the following three dimensions which, although they are not independent, call for separate attention:

- The institutional context (of university vs. school)
- The difference in *the subject's role* within the institution (a student in university or school, vs. a teacher of school mathematics)
- The difference in *mathematical contents* (elementary vs. scientific)

In his book, Klein mainly focuses on the last dimension, and the solutions he proposes can be described as "building bridges" at the level of contents (sometimes with explicit advice to teachers on how to expose a subject). Some of Klein's general proposals have been realized in the course of the twentieth century. For instance, functions and calculus have become a central part of upper secondary education in most European countries. Klein only considers the institutional dimension in the introductory remarks, where he does not fail to point out some basic and problematic discrepancies in basic aims of the two institutions:

> For a long time … university men were concerned exclusively with their sciences, without giving a thought to the needs of the schools, without even caring to establish a connection with school mathematics. *(ibid., p. 1).*

One can safely say that this problem, which Klein describes in past tense in 1908, is not less important today (see e.g. Cuban, 1999). In fact, the causes for university mathematicians' lack of concern and contacts with school mathematics have increased as a consequence of the evolution of institutions: the workforce and institutional frameworks of mathematical research have expanded tremendously since the days of Klein, and mathematics programs and courses prepare students for a still wider range of professions. At the same time, the national institutions of high school have developed in all Western countries, not least in terms of a much larger public and a much wider scope than being a preparatory school for elite university education.

To model Klein's problem, with the three dimensions just stated, we make use of central tools from the anthropological theory of the didactical (ATD), initiated by Chevallard (1991, 1999). Our main justification for this choice of theoretical framework is that it provides us with tools for analysing human practices and institutions which are both applicable and applied by researchers from a wide range of national and scientific backgrounds (as evidenced, for instance, by a series of international congresses on ATD based research, cf. Bosch et al., 2011). Also, it is used in several studies of university mathematics education (see Auteur, 2012 for one survey).

In ATD, human knowledge and practice are modelled as praxeologies, which are collections of practice blocks and theory blocks. A practice block consists of a type of tasks and a corresponding technique, which can be used to accomplish the tasks of the given type. A theory block is attached to a family of practice blocks and consists of technology, i.e. discourse about the techniques, and theory, which justifies, explains and unifies one or (typically) more technologies. Readers unfamiliar with these notions are invited to consult Barbé et al. (2005, sec. 2) or similar references where they are exposed in more detail. We just emphasize here that while techniques and types of task correspond to each other one-to-one, a theory block serves to explain, unify and justify a larger collection of techniques.

An institution is, roughly speaking, a collection of people who share a collection of praxeologies. An example could be everyone involved in the Danish high school, with its repertoire of teaching and study practices in a variety of disciplines and modes. The example also shows why we said "roughly speaking": an institution certainly involves a number of concrete people, with a variety of positions or roles relative to the praxeologies (in the example, being a student or a teacher); but these people come and go over time, and so it is more correct to say that it is the main types of positions which make up the institutions. The praxeologies of an institution also evolve over time, but the positions of its members (e.g., teachers and students) typically remain so stable that we can nevertheless continue to speak of the same institution. We notice also that people in the same position (e.g., students) may certainly develop somewhat different relations to the praxeologies in which they take part, as a function of their position in the institution; still, we may wish to identify and study a smaller number of typical relations, as does Klein in the famous quote given in section 1.

The notation $R_I(x,o)$ was introduced by Chevallard (1991) to indicate the relation of a position $x$ (within an institution $I$) to a praxeology $o$. The notation is just a compact abbreviation, not a mathematical formula. It enables us to represent, in a compact way, the three dimensions of the "double discontinuity" of Klein, as the two passages (each indicated by an arrow):

$$R_{HS}(s,o) \to R_U(\sigma,\omega) \to R_{HS}(t,o)$$

where: o indicates a mathematical praxeology worked on in high school (*HS*) by teachers (*t*) and students (*s*), while ω indicates a mathematical praxeology, which the students (σ) encounter at university (*U*).

In particular, the second part of Klein's problem consists in the lack of (perceived) relevance of $R_U(\sigma,\omega)$ to $R_{HS}(t,o)$, even for the case where o is similar or perhaps a part of ω. While there are certainly many mathematical praxeologies which the student has to relate to at university, but are not close to anything taught in high school, most of the mathematical praxeologies taught in high school (e.g. the practices and knowledge related to derivatives) find some counterpart at the university, perhaps more theoretical, general etc. (for instance, differentiation is also considered for several variables, complex functions etc.). To identify and exploit these counterparts is one strategy to smoothen the second part of the transition.

In this paper, the term 'capstone course' is used to indicate a study unit which is located towards the end of an academic study program, with the aim of concluding or 'crowning' the experience, and to link academic competence and training with the needs of a professional occupation (cf. Durel, 1993; Winsor, 2009). We note here that while the *term* is mainly common in North America, the *phenomenon* exists, with variations, in other parts of the world as well.

In terms of the above model, a capstone course for future teachers aims to develop relations of type $R_U(\sigma,o)$ while drawing on $R_U(\sigma,\omega)$, and in view of the needs for a future relation of type $R_{HS}(t,o)$. As a capstone course takes place within the university programme, the school as an institution remains distant; but it is clear that the motivation of the course is to achieve relations to school mathematics, which can be useful once the student becomes a teacher. Sometimes, the success and outcome of this endeavour may differ according to the qualities of $R_U(\sigma,\omega)$, such as degrees of autonomy of the student with respect to the techniques, technologies and theories of ω. In some cases, we will notice that a further development $R_U^*(\sigma,\omega)$ of $R_U(\sigma,\omega)$ is needed or at least advantageous to achieve a satisfactory result $R_U(\sigma,o)$; in such a case, the complete course of the students may be described as

$$R_U(\sigma,\omega) \to R_U^*(\sigma,\omega) \to R_U(\sigma,o).$$

Now, our overall research questions – which we will elaborate and treat through case studies – can be formulated as follows :

RQ1. What kinds of (new) relations $R_U(\sigma,o)$ is it useful to build between university students and mathematical praxeologies from secondary school, within a university program in mathematics? – here "useful" is to be understood as arguably useful for the aim of preparing students for secondary school teaching.

RQ2. What are the main obstacles to building these relations, in terms of the relations $R_U(\sigma,\omega)$ to mathematical praxeologies which student already have?

RQ3. What further developments $R_U^*(\sigma,\omega)$ of the students' relations to university mathematical praxeologies may be desirable or necessary in order to achieve the goals identified as answers to RQ1?

These questions are open to interpretation in terms of the broadness of answers to be looked for. A relatively modest interpretation of RQ1, which we shall adopt here, will express "kinds of relation" as students' capability to solve tasks coming from *o* itself or from immediately related mathematical praxeologies ω, typically with ω-tasks being about developing theory blocks of *o*. Then the usefulness of students' capability of solving these tasks will be argued through their direct relevance to solve didactic tasks in high school (HS). With this understanding of RQ1, we shall then seek answers to RQ2 and RQ3 related to concrete tasks which students in a capstone course encounter some or many difficulties with, in terms of their relations to university mathematics which either hinder or could further their chances of solving the types of tasks identified in RQ1. In fact, the answers to RQ2 and RQ3 should then identify concrete challenges for university programs in terms of insufficient or desirable relations $R_U(\sigma,\omega)$ and desirable, more advanced relations $R_U^*(\sigma,\omega)$ – where the meaning of "desirable" could relate to the autonomy of students both as regards the technical and theoretical levels of ω. In fact, it is a common experience that university students often relate to "advanced" praxeologies ω in less than advanced ways, that is with a main focus on handling tasks using techniques taught in the course, and with a mostly passive relation to the theory blocks (cf. Auteur, 2008). A main point in our earlier work (Auteur and Auteur, 2007) on didactical engineering for the teaching of real analysis, can be described as constructing new formats of student work which enable the students to develop a more autonomous relation to the theory blocks of advanced mathematical praxeologies.

3. THE BACKGROUND OF OUR CASES

The empirical cases, presented in the next section, come from three years of observation done in the capstone course *UvMat* (an abbreviation for the Danish equivalent of *Mathematics in a Teaching Context*) at the University of Copenhagen.

**Our context.**

UvMat caters to students in their third year of the B.Sc. program in mathematics and aims to help students relate relevant parts of their academic bachelor courses to high school mathematics, in view of a professional life as teachers. To become tenured as such in Denmark the candidate must meet content matter requirements specified by the Ministry of Education. Study programs in mathematics provide courses covering these requirements by offering a comprehensive general program from which the prospective teachers choose appropriately. Being a general program the focus in most courses is on scholarly progress – much in agreement with Klein's vision of how to attain the higher standpoint from which to consider elementary mathematics. In particular, before UvMat, the students had courses on advanced calculus, linear algebra, and real analysis. UvMat is not mandatory, even for prospective high school teachers. It has between 15 and 30 students each year. Mostly participants do only a minor in mathematics (about two years of the bachelor programme) along with a major in another subject. The failure rate at the final exam is relatively high (15-25%, depending on the year).

The overall course goal of UvMat is to enable the students to work autonomously with subjects of high school mathematics from a higher standpoint, implementing that standpoint in a fashion that makes sense in the classroom. Specifically this is spelled out as (1) competencies in solving more demanding problems within high school mathematics, (2) formulating simple as well as challenging and problem-based questions, (3) relating critically to relevant resources (such as text books, technological tools, and websites), and (4) working with applications of mathematics in other subjects.

Klein's approach addresses the transition $R_U(\sigma,\omega) \rightarrow R_U(\sigma,o)$ almost exclusively at the *theoretical* level of the praxeologies, with little emphasis on techniques for problem solving, formulation of tasks etc. The aims of UvMat, outlined above, also contain a strong emphasis on working with *practical blocks* from relevant high school praxeologies. Psychological and pedagogical aspects are of course important for a teacher, but they are not within the scope of UvMat; they are mainly taken up in the induction program, which new teachers enrol in once they are employed at a high school.

**Methodology.**

The main data for our case studies are students' responses to weekly assignments and final exam problems. The written assignments consisted of seven weekly problem sets, six to be answered by groups of 1-4 students and one to be answered individually. The final exam was also individual.

The analysis of student answers was carried out by creating a coding scheme for each item. Each item was subdivided into minimal subtasks and we then identified the techniques each student had employed (or not employed), as a path to identify the exact nature of their difficulties. A detailed example and more explanation of the method are provided in the appendix. In the next section, we present paradigmatic case from this type of analysis.

In addition to this, we have made informal (but very limited) use of the following types of data: (1) email correspondence with some of the students, to support our interpretation of certain problematic answers; (2) a focus group interview of three experienced high school mathematics teachers, to gauge the relevance of the exam problems to authentic teaching of actual secondary mathematics; (3) a voluntary test based on two of the exam questions, given to a voluntary sample of 23 third year students not enrolled in UvMat and analysed with the same method as mentioned above; this will be briefly referred to in our description of main challenges, in view of assessing the extent to which they come from more general shortcomings of students' relations to relevant praxeologies taught in the university programme.

4. PARADIGMATIC CHALLENGES OBSERVED

The research questions presented above can only be treated meaningfully through patient observation and reflexion on cases. We now present a selection of our observations according to the complexity of the mathematical praxeologies which students work with in the course, in order to develop their relation to the praxeologies. This complexity varies from students' capacity to use and explain standard algebraic techniques and technology, to autonomous research and study on a theoretical calculus topic in view of presenting it at a given (high school) level. This leads us to expose five major challenges, which are not just challenges for a single course but also for an entire university program that is supposed to provide its students with adequate mathematical preparation for teaching calculus at high school. We notice here that while many students certainly succeeded with a number of tasks and challenges in the course, we focus here on the most problematic and difficult challenges for us – and for the students – as identified from the analyses described above.

**Challenge 1: autonomous control of algebraic reasoning**

We begin with a challenge which did not occur to us as important during the first years of teaching the course, but which emerged as a serious surprise during the most recent edition (2012). It concerns students' mastery and control of basic algebraic reasoning – that is, not just manipulating symbolic expressions, but using such manipulation in a correct and transparent way. A simple example, to be considered in the following, is to solve an equation with one unknown by hand, with full control of the meaning of all steps.

We speculate that extensive CAS-use by students could be one reason why this is increasingly a challenge even for students 2 or 3 years into a pure mathematics program. While computer algebra systems (CAS) can certainly solve many types of algebraic tasks with great precision and speed, the user of such programs (whether teacher or student) needs to be able to explain and (at least in simple cases) control both input and output in mathematical terms. For these and other reasons we maintain that university students should have a strong and autonomous command of algebraic techniques and technology, including a capacity to develop valid and clear reasoning involving algebraic operations such as the reasoning involved in solving a variety of algebraic equations.

One may also argue that for a teacher at upper secondary school, it is particularly important to be able to formulate algebraic reasoning in a variety of ways and settings, including the logical subtleties involved in solving equations like

$$\sqrt{2x+12} - 1 = x$$

where a mechanical step-by-step rewriting must be supplied with a firm control of the implications between various forms (for the detailed analysis, see the Appendix). A teacher is also be expected to be able to identify and explain the challenges which her students may face with the task, and to describe the challenge in more general terms involving, for instance, mathematical technology and theory related to implications and (solution) sets.

In the 2012 version of our course, we realized from the beginning that the students' relation to algebra and functions was partly insufficient even at the technical level and, to a larger extent, at the level of technology and theory – for instance, about the meaning of the solution to an equation. At the final exam, students were asked to provide a detailed solution to the equation above and then to identify challenges in the task that could be critical for (high school) students. Here, 4 out of 13 students gave a wrong solution to the equation, and were thus certainly unable to answer the second part. Somewhat alarmed by this, we gave the basic task of solving the equation to an informal sample of 23 other students in their last year of bachelor studies; among them we found even higher rates of failure, both at the technical and technological level. This confirms that a significant number of students who take, or could take UvMat, in fact have a relation to algebraic reasoning which is insufficient for developing and explaining the solution of simple equations.

Throughout the course we encountered many other instances of students' inadequate relations to precalculus algebraic techniques and technology, and this certainly constitutes – as a partial answer to RQ2 – a real obstacle to building more advanced relations to praxeologies involving algebraic techniques, which are prerequisites to most calculus praxeologies.

While one could blame this on the high school they come from (thus on the final state of $R_{HS}(s,o)$) the university can clearly not defend to leave students' relation to algebra in the same state. Our experience suggests that teachers of capstone courses must be prepared to detect and work with this kind of clearly inadequate state of $R_U(\sigma,\omega)$, which in fact call for a remedial course rather than a capstone course, and also they must engage in dialogue with teachers of "main" bachelor courses (e.g. in calculus and linear algebra) regarding the appropriate timing for tackling such problems with basic technical and technological capacity.

**Challenge 2: autonomous use of standard calculus techniques**

With the previous challenge in mind, we now approach students' relation to basic techniques from high school calculus, which include calculations and uses of derivatives as in the following exercise, which also appeared in the final exam of our course and where the technical level is only slightly above what is required in Danish high school:

> Assume that the function $f$ is a solution to the differential equation
> $$\frac{dy}{dx} = \exp(y^3).$$
> a) Show that $f$ is strictly increasing.
> b) Show that $f$ is twice differentiable and compute $f''$.

The main challenge in a) is to notice and use that the derivative of $f$ is strictly positive. For many of our students, an obstacle appears from their previous experience with differential equations: they are used to *solve* the equation before (perhaps) considering the solutions. In this case, their standard technique (separation of variables) leads to an integral which cannot be computed in closed form – indeed, a few boldly take this route and get stuck, while almost half of the students don't answer or state they "can't solve it" or the like. As a result, about half of students did not produce a correct answer to this question. The "solving reflex" is clearly counterproductive.

Question b) requires a somewhat challenging use of the chain rule on $\exp(y^3)$, where the challenge is that this expression should be derived with respect to a free variable in $y = y(x)$ which is implicit. At the exam, 22% of the students were able to do this correctly, and in the test with a control group of students at the same level (but who did not attend the course) it was only 17%. These low numbers may also reflect, at least in part, the "solving reflex" obstacle discussed above.

**Challenge 3: autonomous use of instrumented calculus techniques**

In Danish high school, calculus currently involves massive use of computer algebra systems (CAS) such as *Maple* or *TI-Nspire*. In terms of praxeologies, these devices offer techniques (called *instrumented techniques*) that allow the user to solve tasks such as equations, computing limits, plotting graphs and so on. This means that high school praxeologies include instrumented techniques not just as options, but also as students' preferred or (more rarely) unique techniques for many types of tasks. In the mathematics program of the University of Copenhagen (*U*), calculus praxeologies ω involve the use of *Maple* (a professional CAS) in the first semester, but the main goal for $R_U(\sigma,\omega)$ in the program is on developing a closer and more precise relation to the theoretical level (cf. also Auteur, 2011). This includes, for instance, appropriate use of precise definitions, providing details or explanations of proofs etc. As a result, instrumented techniques are much less dominant in students' calculus practices at university.

In the course, we revisit high school level praxeologies to explore the effects and potentials of instrumented techniques, while making use of the theoretical knowledge obtained at *U*. The transition $R_U(\sigma,\omega) \rightarrow R_U(\sigma,o)$ aimed for involves at least two parts: (1) a more subtle use of instrumented techniques along with non-instrumented ones in order to explore difficult tasks and discuss their theoretical perspectives; and (2) to develop the comprehensive and explicit knowledge of instrumented techniques which is needed to design tasks within an instrumented learning environment, and to explain and assess results.

The first part is most directly based on $R_U(\sigma,\omega)$ and we will concentrate on it here. The more delicate relation to instrumented techniques, which we aim for, is better learned at university. In particular, students should learn to combine instrumented techniques with non-instrumented ones, and achieve a better balance between the following aspects of CAS-use:

- The *technical use*, where the instrumented technique is just an easy way to get certain tasks done (this is main the role of CAS to high school students, who sometimes view all tasks of a given type as *either* "very easy", when the instrumented technique works, *or* "impossible", when it doesn't work);
- A *technological use*, to explain and present results using CAS (for instance, to produce illustrative graphs or tables; this occurs more rarely in high school and university, although pc-based CAS-tools have increased the ease of integration of CAS output with normal text);
- A *theoretical use*, such as using CAS as an experimental tool, to investigate a more abstract problem, typically with instrumented techniques being used as complements to pen-and-paper techniques (this occurs in university albeit rarely, see for instance Auteur, 2011).

This more balanced use of CAS is an important example of the new relations $R_U^*(\sigma,\omega)$ at which we aim, especially in the context of praxeologies ω in which students have little or no experience with instrumented techniques.

*Investigating rational functions*

As an example of the difficulties this aim meet with, we consider an item from the weekly assignments :

> Maple gives $a = 1.414213562$ as the 10-digits decimal expansion of $\sqrt{2}$. Investigate the functions $f(x) = \dfrac{x^2 - 2}{x - a}$ and $g(x) = \dfrac{x^2 - 2}{x - \sqrt{2}}$ numerically, algebraically and geometrically. Explain essential differences between the two functions.

A university praxeology ω generated by this task may involve techniques and technology related to poles, removable singularities, polynomial division, and density of ℚ in ℝ. The analogous high school praxeology *o* involves techniques and technology related to vertical asymptotes, zooming in on plots, and numerical tables. The first part of the item aims to supply a basis for relating the two.

In the geometric investigation most students used *Maple* to give plots showing that the graph of *f* has a vertical asymptote whereas the graph of *g* is linear. Some students demonstrated the need to zoom in, by giving plots where "graphs appear identical" and plots where "graphs are clearly different" (students' wording), and one group noted that the necessary degree of zooming is related to the accuracy of the decimal expansion. This was mostly satisfactory.

All students did the numerical investigations using tables of function values. They are generally more unfocused than the geometric investigations or even off the point (for instance using the values $x = 1, 2, ...,10$). Several students interpret "numerical" plainly as a table in itself, mixing exact and floating-point numbers. One student comments on his table: "it is difficult to see that *f* is unbounded". But there are also examples of tables that clearly show this.

Instrumented techniques were not used in students' algebraic analysis of the functions, and while that analysis was adequate, the students missed opportunities (not least for future teachers) related to algebraic CAS-techniques and to coordination with the geometric and numeric investigations based on CAS.

*Modelling with differential equations*

Differential equations represent another calculus topic where we meet the potential and need for all three aspects of CAS use. In one UvMat assignment, we focused on autonomous differential equations, exemplified by fish catch models

$$\frac{dN}{dt} = kN(K - N) - F \qquad \text{(FCM)}$$

where *N* denotes the population size at time *t*, the constants *k* and *K* are model parameters to be interpreted, and *F* is catch per time unit. Investigating such models qualitatively (using CAS) clearly requires one to go beyond the technical use in which one seeks solutions in closed form; in fact, the infinity of closed form solutions to (FCM) may generally say very little about the properties of the model. Instead, one can use CAS to explore (FCM) by producing *phase diagrams* and *direction fields* as well as concrete *solution curves*, and then discuss the relation between the model (FCM) and these three types of diagram. The role of instrumented techniques consequently differs from the example above, as they become essential to the theory blocks, and will be so also in a *HS* learning environment. We denote by *o* the corresponding praxeology.

To produce and make use of the three types of diagram requires new instrumented techniques. These were introduced in a lecture based on an interactive *Maple* sheet with integrated mathematical text, recalling also basic knowledge about differential equations and more basic instrumented techniques such as numerical and symbolic solution commands. In the subsequent project assignment, the students worked with two special cases of the equation (FCM).

The target relation $R_U(\sigma, o)$ requires that the three types of diagrams are *interpreted* and *related* to one another, as well as to the differential equation. The general challenge concerning students' productive of coherent, reasoning technology (already present in Challenge 1) is accentuated through the use of instrumented techniques. For instance, a relatively high-performing student had produced the plot of solution curves on top of a direction field (Fig. 1) to illustrate typical behaviour with respect to equilibrium states. This requires appropriate choices of initial conditions and selections of the plot dimensions.

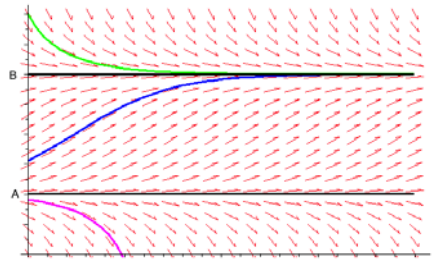

Figure 1. Maple illustration produced by student (cf. text).

Later in the text the student makes a typing error in *Maple* and therefore gets a wrong form ($N(t) = 1 + ce^{-2t}$) of the complete solution; but he does not notice that this solution set is qualitatively very different from the set of solution curves in Fig. 1. The challenge of coordination and integration between CAS-output and non-instrumented techniques, is very visible here and in similar student errors.

**Challenge 4: autonomous work with calculus theory**

The main topics in calculus – such as limits, derivatives and integrals – involve two groups of praxeologies: algebraic praxeologies (based on methods and rules for calculating limits, derivatives etc.) and topological praxeologies (based on existence problems and definitions of limits, derivatives etc.). It is a common trait of secondary level calculus to treat almost exclusively the algebraic praxeologies, with heavy use of instrumented techniques. From the point of view of academic mathematics, this means that student work with "finding" objects (such as limits) for which they have no formal definition or criteria of existence. We refer to Barbé et al. (2005) and Auteur (to appear) for a more detailed discussion of this point.

*A problem on integrals*

Integration is perhaps the most advanced topic that is dealt with in high school mathematics. The praxeologies *o* taught there involve informal explanations of what the definite integral computes (certain areas, averages etc.) and how to compute it (antiderivatives, instrumented techniques); at the most advanced levels, we also find sketchy arguments of the fundamental theorem of calculus, bypassing any serious criteria for the existence of the integral. At the university, the students are given rigorous definitions of integrals (in fact, various alternatives), as part of a more or less comprehensive treatment of real analysis. The question, naturally, is what students σ retain from these presentations and if they are able to make use of $R_U(\sigma,\omega)$ to build a new relation $R_U(\sigma,o)$ which is relevant to the teaching of *o* in high schools. A minimal interpretation of Klein's ideal of teaching "from a higher standpoint" is that students can reason autonomously about integrals in a mathematically sound way, so that they will not simply *forced* to resort to unreasoned statements once they become teachers. For instance, they should be able to explain how the Riemann integral gives sense to the area of certain subsets of the plane, and they should be able to reason autonomously and correctly about basic notions in *o* such as integral, integrability, antiderivative and continuity. As we shall see, this is far from guaranteed from the fact that students have been *presented* to all of this (and much more!) in undergraduate courses.

Here is an example of a task given to the students in UvMat:

Let the function *f* be integrable on the interval *I* = [0,1], and define the function *F* by

$$F(x) = \int_0^x f(t)dt, \quad x \in I.$$

> Justify with rigorous reasoning whether or not it follows from what is given that $F$ is continuous.
>
> Justify, likewise with rigorous reasoning, whether or not one can deduce that $F$ is differentiable.

The questions are open ended, asking for not just a proof, but also an answer to be proved. None of the students gave a fully satisfactory solution, while about 1/3 came close. Two types of shortcomings were prevalent among the rest. One is to present a sequence of statements which involve relevant notions and locally appear sound, but with little logical coherence among the statements. The following gives an impression of what that kind of "reasoning" may look like (we stress that we copied and translated the student production *exactly*, including what appears as accidental errors):

> ... Let $\varepsilon > 0$. We look at: Let $x_0 \in I. \mid F(x_0) - F(x) \mid = \int_x^{x_0} f(t)dt$, by the subdivision rule. $\int_x^{x_0} f(t)dt \to 0$ as $x_0 \to 0$. Choose $\delta$ so that $\mid x_0 - x \mid < \delta \Rightarrow \int_x^{x_0} f(t)dt < \varepsilon$
>
> so that it is concluded that $F$ is continuous ...

The other shortcoming is the widespread mistake to consider the integral $F(x) = \int_c^x f(t)dt$ to be *defined* as the value at $x$ of the *antiderivative* of $f$, which satisfies $F(c) = 0$ (this turns the Fundamental Theorem of Analysis into a circular statement). Students produce disguised versions of this, as in the following student answer to b):

> By the definition of antiderivative we get that if $F(x) = \int_c^x f(t)dt, \quad x \in I,$ then $F$ is differentiable in each point $x \in I$ with derivative $F'(x) = f(x)$.

This kind of solution both reflects a resisting "definition" (or explanation) from high school, where the only explicit formula for the integral is the one involving an antiderivative of the integrand, and an insufficient experience of students with autonomous use of the notions of integrability and continuity. Such basic shortcomings are difficult to deal with in a capstone course where the higher viewpoint is, at least to some extent, assumed.

*Construction of the exponential function*

An obvious and potentially rich topic in a capstone course like UvMat is to study and indeed define the meaning of exponential expression $a^x$ for $a \geq 0$ and, crucially, with $x$ an arbitrary real number. In high school, such expressions are introduced very early as examples of elementary functions, which will later become central examples and building blocks in the study of calculus. The challenge for the high school teacher is to give some meaning to $a^x$ in the absence of any rigorous theory of real numbers and previous to the study of limits and other elements of the calculus. The common approach is to give a more or less detailed algebraic justification of the formula $a^{m/n} = \sqrt[n]{a^m}$, and then merely claim that one can extend this definition from rational to arbitrary real numbers. Here are some explanations of this extension in Danish high school textbooks:

> The power is calculated by approximating the exponent by a finite decimal number. How many decimals you include depend on the required accuracy *(Timm & Svendsen, 2005, 26; translated from Danish by the authors)*

> In Chapter 3 we saw how to calculate powers where the exponent is integer and positive, 0,

integer and negative, and rational (fraction). Strictly speaking we have not explained the meaning of a symbol like $7^{\sqrt{3}}$ but we assume CAS will take care of this. *(Carstensen, Frandsen & Studsgaard, 2006, 82; translated from Danish by the authors)*

For the following discussion, we denote by ω the mathematical praxeology built from the university experience, in view of establishing exponential functions $a^x$ and their basic properties, with a complete mathematical theory to justify the extension to real numbers. In the high school version *o* of ω, the exponential functions will certainly have to be introduced on the domain of real numbers; but at the moment where they are introduced, the available theory excludes a rigorous theoretical justification, so that any justification will need to be somewhat informal, as the above examples suggest. Still, the quality of explanation and activities proposed by teachers could certainly vary a lot, from the worst (not recognizing the problem or believing it to be insignificant) to the better (proposing a range of activities and explanations which could even anticipate or prepare a more rigorous work with approximations and limits of functions). Our experience shows that university students' relation to the theoretical level of ω is very weak: when informally asked, it appears that most of our students have never realized that the general definition of $a^x$ poses a problem, and none have ever seen a complete solution. Here, different "complete solutions" could be more or less helpful for the future teacher. For instance, we believe that Klein's coverage of the theme (Klein, 1908, pp. 144-162), in terms of logarithms and complex functions, misses the subtleties they will be able to address, even if partially, in high school. Instead, we propose to our students a study of the following key lemma (with free use of theoretical results they have encountered in real analysis):

**Lemma.** If *a* is a positive real number, *x* is a real number and $(r_n)$ is a sequence of rational numbers converging to *x*, then the sequence $(a^{r_n})$ is convergent. Moreover if $(q_n)$ is another sequence of rational numbers converging to x, then the limit of the sequence $(a^{q_n})$ is the same as the limit of $(a^{r_n})$.

Any proof will have to invoke deep properties of the real numbers, in particular completeness in some form. (For one proof, see Bremigan, Bremigan & Lorch, 2011, 294-295).

One of the challenges for students in studying such proofs lie in the fact that the knowledge which they are supposed to use are not always what is most familiar to them. For instance, one task given to students was to prove, through a number of steps, the following special case of the above lemma: *If a is a positive real number, and if $(r_n)$ is any sequence of rational numbers converging to 0, then the sequence $(a^{r_n})$ converges to* 1. Several students at some point invoked properties of exponential functions, such as continuity (then, the proof is of course trivial); instead they should use properties of the exponential function defined on the set of rational numbers, along with general properties of the real numbers. As teachers we had not foreseen the difficulty or need of making explicit what theoretical elements it makes sense to use when the task is to construct the function from scratch. On the other hand, it may not surprise that some students found it weird that they *could and should* use (for them) more advanced theoretical knowledge related to order and convergence in ℝ, but were criticized for using familiar results like the continuity of exponential functions. The students' lack of autonomy with the theoretical level of ω had led us to construct theoretical tasks in which students had to carry out smaller steps, rather than to construct and explain the wider coherence and purpose. The main conclusion of the challenge which even this represents is that *students need to be given more and earlier experiences at university with autonomous construction of theoretical structures, involving definitions, partial results, proofs etc.* To be manageable, this would preferably be done in the setting of mathematical praxeologies where students are already very familiar with the practice blocks, as can be expected in the case of ω considered here; also, one will need a reasonable progression in the size of those praxeologies, and considerable time to share reflections and validate their ideas.

**Challenge 5: autonomous study beyond text books**

Autonomous search and use of mathematical text can be seen as one of the most advanced and difficult relationship of type $R_U(\sigma,o)$ which one may seek to develop in a capstone course, based on $R_U(\sigma,\omega)$ where ω is one or more mathematical organisations from the university curriculum, applicable to the study of a high school mathematical organisation $o$. To exemplify this, consider the technique of least squares for simple linear regression – regularly taught in high school. The type of task generating $o$ consists in determining, for a data set $\{(x_1,y_1),...,(x_n,y_n)\} \subseteq \mathbb{R}^2$, the line $y = ax + b$, which approximates the data points best in the sense that the square sum $S(a,b) = \sum_{k=1}^{n}(ax_k + b - y_k)^2$ is minimized. The technique is just a formula (often implicit, with an instrumented technique); but how to justify it in a way accessible to high school students?

Notice first that relevant university level theory block could come from both calculus and linear algebra parts of the undergraduate curriculum. Many university textbooks use calculus, and one easily finds that $S$ has a unique critical point; it takes more technicalities to actually prove that this point is a minimum. A similar application of linear algebra can be obtained using the vector projection of $\vec{y}$ onto the two-dimensional subspace of $\mathbb{R}^n$ spanned by $\vec{x} = (x_1,...,x_n)$ and $\vec{1} = (1,...,1)$; the solution is simply the coordinates of that projection in the basis $\{\vec{x},\vec{1}\}$. But high school theory blocks do not include partial derivative tests or vector projections. Alternatives are not found in Danish high school text books, and they are not (as of early 2013) so easily found on the Internet; but a slightly more insisting Google search does lead one to texts like Key (2005), with an elementary proof based on "completion of squares".

During the 2011 of the course, we asked students to autonomously search for a proof, which is within the reach of high school mathematics (i.e. common theory blocks of a suitable $o$). This was impossible for all students, even with some help for getting started on the search; instead they all attempted to elementarise the calculus proof as an "analogue" of a one variable optimization problem.

In fact, students get very little experience and aptitude from their undergraduate studies when it comes to autonomous search and study of mathematical literature, in view of solving a concrete problem (such as finding an alternative proof, getting ideas for posing exercises etc.). Finding ways to change this aspect of students' relation to mathematical praxeologies in this direction appears to us a main open problem in teaching this course. The problem is not only that they have relatively little experience with autonomous search for resources, especially with constraints like "find a solution accessible for high school students" (as in the case above). The difficulty lies also, it seems, in distinguishing potentially useful resources from irrelevant ones, and especially in working with the first type of resources as one rarely find a complete "solution" in one resource. In the case above, even after finding a text like Key (2005), details need to be worked out by the students to realize whether or not one has found what one was looking for – and then to work out the details to give an explicit and personalized solution, as commonly required in mathematics courses. It should be mentioned that in 2012, we simply gave the text by Key to students, and most of them were then able to accomplish that last step.

Based on this and similar experience, we believe that a development of students' relation to mathematical praxeologies to include capacities for autonomous search and use of mathematical literature will be a difficult challenge for any capstone course departing from standard, text book based bachelor programs in pure mathematics. We also think that it is one of the most important aims in capstone courses for future teachers, given that teachers can – sometimes must – work with a wide variety of resources (cf. Gueudet & Trouche, 2009) – which nowadays naturally include Internet based ones. This aim is particularly evident in view of a context like Danish high school where significant parts of the teaching are done as supervision of individualized and multidisciplinary student "projects".

5. CONCLUSIONS

In this paper, we have uncovered a number of requirements and obstacles to construct a "bridge" between the mathematical praxeologies of contemporary university and high school, in the context of a capstone course where high school level calculus is studied and put in perspective using technical and theoretical elements of the university program. Certainly the basic "double discontinuity" identified by Klein remains for the students whose academic preparation for high school teaching is mainly based on university mathematics studies. In fact, the institutional and societal conditions have clearly changed considerably in 100 years, in ways that tend to widen the gap to be bridged, as explained in Sec. 2.

We do share Klein's belief that the advanced standpoint of university mathematics is potentially relevant to develop teachers' "deep" knowledge of the high school subject. But most of our experiences confirm that the challenges we face today are acute, both in terms of the shortcomings of students' grasp of relevant university mathematics, and in terms of the difficulty of creating situations where they can experience and realize the relevance of what they do know in solving problems related to calculus as taught in high school.

As regards calculus praxeologies taught in high school, the increasing importance of instrumented techniques, as well as the informal nature of theoretical blocks, requires specific attention to the future high school teacher's preparation at university: he must be prepared to solve and construct tasks for his students which go beyond a sequence of unrelated, meaningless procedures which in the end amount to choosing relevant commands or buttons of a CAS device. He must, in particular, know a number of alternative approaches to the "hard" topological problems that, in calculus, are based on the completeness of real numbers - such as the definition and existence of central objects like elementary classes of functions, limits, derivatives and integrals. Some of these approaches should make use of the potential of CAS to visualize and compute, while others must be based on simplified heuristic arguments and shortcuts which do not simply amount to circularity, mysterious beliefs or story telling.

At the same time, and in a sense before that, we must also face – and further investigate – a number of shortcomings in some university students' relation to mathematical practices which are or ought to be central also in their previous studies, such as the capacity to solve equations correctly (challenge 1), to make use of calculus techniques in non-standard problems (challenge 2), to make use of instrumented techniques to investigate concrete problems (challenge 3), to device simple but correct mathematical arguments (challenge 4), and to search and study mathematical literature in an autonomous way, for instance to identify alternative proofs under boundary conditions for the machinery to be deployed (challenge 5). In fact, all of these capacities should be established in analysis and algebra courses before "visits" to advanced topics more distant to high school calculus. And even when that is done, so that "remedial measures" become less needed, Klein's problem remains and motivates a continuing effort to develop the contents and methods of capstone courses.

APPENDIX : DETAILED EXAMPLE OF ITEM ANALYSIS

The first question in the exam item discussed in "Challenge 1" is the following :

*Solve the equation* $\sqrt{2x+12} - 1 = x$ ; *provide all intermediate steps of your solution.*

Our coding of answers is based on the degree to which seven subtasks (given below) were identified and solved, and also on the explicit connections provided between them and the given task. We stress that there is no contention that these subtasks should be considered in a particular order or that all of them need to be addressed:

Subtask 1: Decide for what values of *x* the equation makes sense (namely $x \geq -6$)

Subtask 2: Rewrite the equation as $\sqrt{2x+12} = x+1$

Subtask 3: Infer that the result from task 2 *implies* that $x \geq -1$,

Subtask 4: Infer that the result from subtask 2 *implies* that $2x+12 = (x+1)^2$

Subtask 5: Observe that $x \geq -1$ and $2x+12 = (x+1)^2$ *together* imply $\sqrt{2x+12} = x+1$ (noting that $x \geq -1$ also ensures the condition from subtask 1) so that the given equation is logically *equivalent* to $2x+12 = (x+1)^2 \wedge x \geq -1$

Subtask 6: Solve $2x+12 = (x+1)^2$ (solutions: $x = \pm\sqrt{11}$)

Subtask 7: Identify the complete solution to the equation as $x = \sqrt{11}$.

We then created a table with one row for each respondent and one column for each subtask; the presence and character of solutions to the subtasks was then indicated, for each student, using a specific coding system.

It is important to note that our analysis of items was not "fixed" but could be changed to admit alternative, correct solutions; but for almost all items, we did manage to predict the steps taken by students. In the present case, this does not mean that they solved all subtasks or did them in the given order. In fact, a common technique was to solve subtask 2, 4 and 6 in that order, and then either (incorrectly) state $x = \pm\sqrt{11}$ as the solution, or insert the original equation and observe that only $x = \sqrt{11}$ "works", so that this is the solution. This final step was coded as a partial solution of subtask 3 as it leads to observe that one concrete number less that –1 does not satisfy the equation; with the other steps, and of course an explicit argument, this is indeed a valid solution.

Clearly, the assessment of the task is not completed by using this coding, which amounts merely to identifying the students' *technique* ; their *technology* and *theory* is then identified through explicit connections between the steps (as indicated in the main text on Challenge 1), explicit appeal to specific results, erroneous inferences, etc. On the other hand, the table resulting from the coding related to subtasks was indeed very useful to provide an overview of students' capacities and challenges related to a given task.